# CONCENTRATION AROUND THE MEAN FOR MAXIMA OF EMPIRICAL PROCESSES

By T. Klein and E. Rio

*Université de Versailles Saint Quentin en Yvelines*

In this paper we give optimal constants in Talagrand's concentration inequalities for maxima of empirical processes associated to independent and eventually nonidentically distributed random variables. Our approach is based on the entropy method introduced by Ledoux.

**1. Introduction.** Let $X_1, X_2, \ldots$ be a sequence of independent random variables with values in some Polish space $\mathcal{X}$ and let $\mathcal{S}$ be a countable class of measurable functions from $\mathcal{X}$ into $[-1, 1]^n$. For $s = (s^1, \ldots, s^n)$ in $\mathcal{S}$, we set

$$S_n(s) = s^1(X_1) + \cdots + s^n(X_n). \tag{1.1}$$

In this paper we are interested in concentration inequalities for $Z = \sup\{S_n(s) : s \in \mathcal{S}\}$.

Now let us recall the main results in this direction. Starting from concentration inequalities for product measures, Talagrand (1996) obtained Bennett type upper bounds on the Laplace transform of $Z$ via concentration inequalities for product measures. More precisely he proved

$$\log \mathbb{E} \exp(tZ) \leq t\mathbb{E}(Z) + Vab^{-2}(e^{bt} - bt - 1) \tag{1.2}$$

for any positive $t$. Here

$$V = \mathbb{E}\left(\sup_{s \in \mathcal{S}} \sum_{k=1}^n (s^k(X_k))^2\right).$$

In order to analyze the variance factor $V$, set

$$V_n = \sup_{s \in \mathcal{S}} \operatorname{Var} S_n(s). \tag{1.3}$$

---

Then, one can derive from the comparison inequalities in Ledoux and Talagrand (1991) that $V_n \leq V \leq V_n + 16\mathbb{E}(Z)$ [see Massart (2000), page 882]. Consequently $V$ is often close to the maximal variance $V_n$. The conjecture concerning the constants is then $a = b = 1$. The constant $a$ plays a fundamental role; in particular, for Donsker classes, $a = 1$ gives the exact rate function in the moderate deviations bandwidth. Nevertheless it seems difficult to reach $a = 1$ via Talagrand's method [see Panchenko (2001) for more about the constants in Talagrand's concentration inequalities for product measures]. In order to obtain concentration inequalities more directly, Ledoux (1996) used a log-Sobolev type method together with a powerful argument of tensorization of the entropy. When applied to $\exp(tZ)$, this method yields a differential inequality (this is the so-called Herbst argument) on the Laplace transform of $Z$ and gives (1.2) again. Applying Ledoux's method, Massart (2000) obtained $a = 8$ in (1.2) with Talagrand's variance factor and $a = 4$ in (1.2) with the variance factor $V_n + 16\mathbb{E}(Z)$. Later on, Rio (2002) proved (1.2) for independent and identically distributed (i.i.d.) random variables (in the i.i.d. case $s^1 = \cdots = s^n$) with $a = 1$, $b = 3/2$ and a variance factor $v = V_n + 2\mathbb{E}(Z)$. Next, Bousquet (2003) found a nice trick to improve Rio's inequality. He proved (1.2) with $a = b = 1$ and the variance factor $v$ in the i.i.d. case. For negative values of $t$, Klein (2002) obtained (1.2) in the i.i.d. case with $a = 1$, $b = 4$ and the same factor $v$.

Here we are interested in optimal constants in Talagrand's inequalities for nonidentically distributed random variables. Our approach to obtain the best constants is to apply the lemma of tensorization of the entropy proposed by Ledoux (1996). However, the differential inequality on the Laplace transform of $Z$ is more involved than in the i.i.d. case. Therefore the results are suboptimal in the large deviations bandwidth. We start by right-hand side deviations.

THEOREM 1.1. *Let $\mathcal{S}$ be a countable class of measurable functions with values in $[-1,1]^n$. Suppose that $\mathbb{E}(s^k(X_k)) = 0$ for any $s = (s^1, \ldots, s^n)$ in $\mathcal{S}$ and any integer $k$ in $[1, n]$. Let $L$ denote the logarithm of the Laplace transform of $Z$. Then, for any positive $t$,*

(a) $$L(t) \leq t\mathbb{E}(Z) + \frac{t}{2}(2\mathbb{E}(Z) + V_n)(\exp((e^{2t} - 1)/2) - 1).$$

*Consequently, setting $v = 2\mathbb{E}(Z) + V_n$, for any positive $x$,*

(b) $$\mathbb{P}(Z \geq \mathbb{E}(Z) + x) \leq \exp\left(-\frac{x}{4}\log(1 + 2\log(1 + x/v))\right)$$

*and*

(c) $$\mathbb{P}(Z \geq \mathbb{E}(Z) + x)$$
$$\leq \exp\left(-\frac{x^2}{v + \sqrt{v^2 + 3vx} + (3x/2)}\right) \leq \exp\left(-\frac{x^2}{2v + 3x}\right).$$

REMARK 1.1. In the spirit of Massart's paper (2000), Theorem 1.1(b) can be improved for large values of $x$ to get a Bennett type inequality with $a = 1$.

REMARK 1.2. Theorem 1.1 applies to set-indexed empirical processes associated to nonidentically distributed random variables. In that case $s^i(X_i) = \mathbb{1}_{X_i \in S} - \mathbb{P}(X_i \in S)$ and consequently the centering constant depends on $i$. Some different concentration inequalities for set-indexed empirical processes are given in Rio [(2001), Theorem 4.2 and Remark 4.1]. However, due to the concavity of the polynomial function $u(1-u)$, the variance factor in Rio (2001) is suboptimal for nonidentically distributed random variables. Here, as a by-product of Theorem 1.1(a), we get the upper bound below for the variance of $Z$.

COROLLARY 1.1. *Under the assumptions of Theorem* 1.1(a), $\operatorname{Var} Z \leq V_n + 2\mathbb{E}(Z)$.

For left-hand side deviations, the concentration bounds are similar. However, the proof is more intricate. We emphasize that the proof of Theorem 1.1 is not relevant for left-hand side deviations. This is the reason why we need to compensate the empirical process for left-hand side deviations.

THEOREM 1.2. *Under the assumptions of Theorem* 1.1, *for any positive* $t$,

(a) $$L(-t) \leq -t\mathbb{E}(Z) + \frac{v}{9}(e^{3t} - 3t - 1).$$

*Consequently, for any positive* $x$,

(b) $$\mathbb{P}(Z \leq \mathbb{E}(Z) - x) \leq \exp\left(-\frac{v}{9}h\left(\frac{3x}{v}\right)\right),$$

*where* $h(x) = (1+x)\log(1+x) - x$, *and*

(c) $$\mathbb{P}(Z \leq \mathbb{E}(Z) - x) \leq \exp\left(-\frac{x^2}{v + \sqrt{v^2 + 2vx} + x}\right) \leq \exp\left(-\frac{x^2}{2v + 2x}\right).$$

REMARK 1.3. Theorem 1.2(b) improves on Theorem 1.1, inequality (2) in Klein (2002). However, Klein gives additional results for functions with values in $]-\infty, 1]$ and subexponential tails on the left [cf. inequality (3), Theorem 1.1].

Let us now apply Theorems 1.1 and 1.2 to randomized processes, as defined in Ledoux and Talagrand [(1991), Section 4.3]. Let $X_1, X_2, \ldots, X_n$ be a sequence of independent and centered random variables with values in $[-1,1]$. Let $T$ be some countable set and let $\zeta_1, \zeta_2, \ldots, \zeta_n$ be numerical functions on $T$. Let

$$Z = \sup\{X_1\zeta_1(t) + X_2\zeta_2(t) + \cdots + X_n\zeta_n(t) : t \in T\}.$$

The random variable $Z$ corresponds to the class of functions $\mathcal{S} = \{s_t : t \in T\}$, where the components $s_t^i$ of $s_t$ are defined by $s_t^i(x) = x\zeta_i(t)$. Assuming that

$$V_n = \sup_{t \in T} \sum_{k=1}^n \zeta_k^2(t)\mathbb{E}(X_k^2) < \infty \quad \text{and} \quad M = \sup_{k \in [1,n]} \sup_{t \in T} |\zeta_k(t)| < \infty,$$

Corollary 1.1 gives $\operatorname{Var} Z \leq V_n + 2\mathbb{E}(Z)$. Let us compare this variance bound with the known results. Theorem 3 in Bobkov (1996) applied to $Z$ yields $\operatorname{Var} Z \leq 2V$, where

$$V = \mathbb{E}\left(\sup_{t \in T} \sum_{k=1}^n (\zeta_k(t)X_k)^2\right)$$

is Talagrand's variance factor. If the random variables $X_1, X_2, \ldots, X_n$ are symmetric signs, then $Z$ is the maximum of a Rademacher process and $V = V_n$. In that case Corollary 1.1 improves the known bounds on $\operatorname{Var} Z$ as soon as $2\mathbb{E}(Z) < V_n$. For Rademacher processes, the concentration inequality (4.10) in Ledoux and Talagrand (1991) yields

(1.4) $$\mathbb{P}(Z \geq m_Z + x) \leq \exp(-x^2/(8V_n)),$$

where $m_Z$ denotes a median of $Z$. Theorems 1.1 and 1.2 provide exponential bounds with a factor 2 instead of 8. However, our variance factor is greater than $V_n$ and our bounds are not sub-Gaussian. Finally, we refer the reader to Bousquet (2003) or Panchenko (2003) for concentration inequalities (with suboptimal variance factor) for randomized or empirical processes in the unbounded case.

**2. Tensorization of entropy and related inequalities.** In this section we apply the method of tensorization of the entropy to get an upper bound on the entropy of positive functionals $f$ of independent random variables $X_1, X_2, \ldots, X_n$.

NOTATION 2.1. Let $\mathcal{F}_n$ be the $\sigma$-field generated by $(X_1, \ldots, X_n)$ and let $\mathcal{F}_n^k$ be the $\sigma$-field generated by $(X_1, \ldots, X_{k-1}, X_{k+1}, \ldots, X_n)$. Let $E_n^k$ denote the conditional expectation operator associated to $\mathcal{F}_n^k$.

In this paper, the main tool for proving concentration inequalities is the following consequence of the tensorization inequality in Ledoux (1996).

PROPOSITION 2.1. *Let $f$ be some positive $\mathcal{F}_n$-measurable random variable such that $\mathbb{E}(f \log f) < \infty$ and let $g_1, g_2, \ldots, g_n$ be any sequence of positive and integrable random variables such that $\mathbb{E}(g_i \log g_i) < \infty$. Then*

$$\mathbb{E}(f \log f) - \mathbb{E}(f) \log \mathbb{E}(f)$$
$$\leq \sum_{k=1}^{n} \mathbb{E}(g_k \log(g_k/E_n^k g_k)) + \sum_{k=1}^{n} \mathbb{E}((f - g_k) \log(f/E_n^k f)).$$

PROOF. Set $f_k = E_n^k f$. By the tensorization inequality in Ledoux (1996),

$$(2.1) \qquad \mathbb{E}(f \log f) - \mathbb{E}(f) \log \mathbb{E}(f) \leq \sum_{k=1}^{n} \mathbb{E}(f \log(f/f_k)).$$

Now

$$(2.2) \qquad \mathbb{E}(f \log(f/f_k)) = \mathbb{E}(g_k \log(f/f_k)) + \mathbb{E}((f - g_k) \log(f/f_k)).$$

Since $E_n^k(f/f_k) = 1$, we have

$$\mathbb{E}(g_k \log(f/f_k)) \leq \sup\{\mathbb{E}(g_k h) : h \ \mathcal{F}_n\text{-measurable}, E_n^k(e^h) = 1\}.$$

Hence, from the duality formula for the relative entropy in Ledoux (1996),

$$\mathbb{E}(g_k \log(f/f_k)) \leq \mathbb{E}(g_k \log(g_k/E_n^k g_k)).$$

Together with (2.2), it implies that

$$(2.3) \quad \mathbb{E}(f \log(f/f_k)) \leq \mathbb{E}(g_k \log(g_k/E_n^k g_k)) + \mathbb{E}((f - g_k) \log(f/f_k)),$$

and Proposition 2.1 follows. □

**3. Right-hand side deviations.** To prove Theorems 1.1 and 1.2, we start by proving the results for a finite class of functions. The results in the countable case are derived from the finite case using the Beppo Levi lemma. Consequently, throughout the sequel we may assume that $\mathcal{S} = \{s_1, \ldots, s_m\}$.

As mentioned in the Introduction, the deviation of $Z$ on the right is easier to handle than the deviation on the left. In fact, for positive $t$, the functional $\exp(tZ)$ is an increasing and convex function with respect to the variables $s_i^k(X_k)$. This is not the case for negative values of $t$. Consequently, upper bounds for the Laplace transform of $Z$ via the Herbst–Ledoux method are more difficult to handle for negative values of $t$. In Section 4, we will introduce compensated processes in order to handle the deviation on the left.

DEFINITION 3.1. Let $\tau$ be the first integer such that $Z = S_n(s_\tau)$. Set $f = \exp(tZ)$ and $f_k = E_n^k(f)$. Let $P_n^k$ denote the conditional probability measure conditionally to $\mathcal{F}_n^k$.

Set

$$g_k = \sum_i P_n^k(\tau = i) \exp(tS_n(s_i)). \tag{3.1}$$

Let $F$ denote the Laplace transform of $Z$. From Proposition 2.1,

$$\begin{aligned}tF'(t) - F(t)\log F(t) \\ \leq \sum_{k=1}^n \mathbb{E}(g_k \log(g_k/E_n^k g_k)) + \sum_{k=1}^n \mathbb{E}((f-g_k)\log(f/f_k)).\end{aligned} \tag{3.2}$$

Since $f - g_k \geq 0$, the upper bound on the second term in (3.2) will be derived from Lemma 3.1.

LEMMA 3.1.  *With the notation of Definition 3.1, $\exp(ts_\tau^k(X_k)) \geq (f/f_k) \geq \exp(-2t)$ a.s.*

PROOF.  Let $S_n^k(s) = S_n(s) - s^k(X_k)$. Let $\tau_k$ be the first integer in $[1,m]$ such that

$$S_n^k(s_{\tau_k}) = \sup\{S_n^k(s) : s \in \mathcal{S}\} =: Z_k. \tag{3.3}$$

Clearly $Z_k$ is an $\mathcal{F}_n^k$-measurable random variable and

$$\exp(tZ_k + t) \geq f \geq \exp(tZ_k)\exp(ts_{\tau_k}^k(X_k)). \tag{3.4}$$

Since the stopping time $\tau_k$ is $\mathcal{F}_n^k$-measurable, $E_n^k(s_{\tau_k}^k(X_k)) = 0$ by the centering assumption on the elements of $\mathcal{S}$. It follows that

$$E_n^k f \geq \exp(tZ_k) E_n^k(\exp(ts_{\tau_k}^k(X_k))) \geq \exp(tZ_k) \geq \exp(tS_n^k(s_\tau)). \tag{3.5}$$

Hence $f_k \geq f \exp(-ts_\tau^k(X_k))$, which implies the left-hand side inequality in Lemma 3.1.

We now prove the second inequality in Lemma 3.1. From the left-hand side inequality in (3.4), $\exp(tZ_k + t) \geq E_n^k(f)$. Next, from the right-hand side inequality in (3.4), $\exp(tZ_k) \leq \exp(tZ + t)$. Hence $f_k \leq f\exp(2t)$, which implies the second part of Lemma 3.1. □

From Lemma 3.1 and the facts that $f - g_k \geq 0$ and $ts_\tau^k(X_k) \leq t$ we get that

$$\mathbb{E}((f-g_k)\log(f/f_k)) \leq t\mathbb{E}(f - g_k). \tag{3.6}$$

We now bound up the first term in (3.2). Set

$$h_k = \sum_{i=1}^m P_n^k(\tau = i)\exp(tS_n^k(s_i)). \tag{3.7}$$

CONCENTRATION FOR EMPIRICAL PROCESSES    7

The random variable $h_k$ is positive and $\mathcal{F}_n^k$-measurable. Hence, from the variational definition of the relative entropy [cf. Ledoux (1996), page 68],

$$E_n^k(g_k \log(g_k/E_n^k g_k)) \leq E_n^k(g_k \log(g_k/h_k) - g_k + h_k).$$

Putting this inequality in (3.2) and using (3.6), we get

$$tF' - F \log F$$
$$(3.8) \quad \leq \sum_{k=1}^n \mathbb{E}(g_k \log(g_k/h_k) + (1+t)(h_k - g_k)) + t \sum_{k=1}^n \mathbb{E}(f - h_k).$$

In order to bound up the second term on the right-hand side, we will use Lemma 3.2.

LEMMA 3.2. *Let $(h_k)_{k \leq n}$ be the finite sequence of random variables defined in (3.7). Then*

$$\sum_{k=1}^n \mathbb{E}(f - h_k) \leq e^{2t} F(t) \log F(t).$$

PROOF. Since the random variables $S_n^k(s)$ are $\mathcal{F}_n^k$-measurable,

$$h_k = E_n^k \left( \sum_{i=1}^m \mathbb{1}_{\tau = i} \exp(t S_n^k(s_i)) \right) = E_n^k(\exp(t S_n^k(s_\tau))).$$

It follows that

$$(3.9) \quad \sum_{k=1}^n \mathbb{E}(f - h_k) = \sum_{k=1}^n \mathbb{E}(f(1 - \exp(-t s_\tau^k(X_k)) - e^{2t} t s_\tau^k(X_k))) + t e^{2t} F'(t).$$

Now, from Lemma 3.1, $t s_\tau^k(X_k) \geq \log(f/f_k) \geq -2t$. Since $1 - \exp(-x) - e^{2t}x$ is a nonincreasing function of $x$ on the interval $[-2t, +\infty[$, it follows that

$$\mathbb{E}(f(1 - \exp(-t s_\tau^k(X_k)) - t e^{2t} s_\tau^k(X_k))) \leq \mathbb{E}(f - f_k - e^{2t} f \log(f/f_k)).$$

From the equality $\mathbb{E}(f_k) = \mathbb{E}(f)$, we get that

$$\mathbb{E}(f - f_k - e^{2t} f \log(f/f_k)) = -e^{2t} \mathbb{E}(f \log(f/f_k)).$$

Hence, summing on $k$ and applying (2.1),

$$\sum_{k=1}^n \mathbb{E}(f(1 - \exp(-t s_\tau^k(X_k)) - t e^{2t} s_\tau^k(X_k))) \leq e^{2t}(F \log F - tF'),$$

which, together with (3.9), implies Lemma 3.2. □

Next, we bound up the first term on the right-hand side in (3.8).

DEFINITION 3.2. Let $r(t,x) = x\log x + (1+t)(1-x)$.

With the above definition

$$g_k \log(g_k/h_k) + (1+t)(h_k - g_k) = h_k r(t, g_k/h_k).$$

From the convexity of $r$ with respect to $x$,

$$h_k r(t, g_k/h_k) \leq \sum_i P_n^k(\tau = i) \exp(tS_n^k(s_i)) r(t, \exp(ts_i^k(X_k))),$$

which ensures that

(3.10) $\quad E_n^k(h_k r(t, g_k/h_k)) \leq \sum_i P_n^k(\tau = i) \exp(tS_n^k(s_i)) \mathbb{E}(r(t, \exp(ts_i^k(X_k)))).$

Here we need the bound below.

LEMMA 3.3. *Let $r$ be the function defined in Definition 3.2. For any function $s$ in $\mathcal{S}$ and any positive $t$,*

$$\mathbb{E}r(t, \exp(ts^k(X_k))) \leq \frac{t^2}{2}\mathbb{E}(s^k(X_k))^2.$$

PROOF. Let $\eta(x) = r(t, e^{tx}) = txe^{tx} + (t+1)(1 - e^{tx})$. We will prove that, for any $x \leq 1$,

(3.11) $\qquad\qquad\qquad\qquad \eta(x) \leq x\eta'(0) + (tx)^2/2.$

Set $\delta(x) = \eta(x) - x\eta'(0) - (tx)^2/2$. Then $\delta(0) = 0$ and $\delta'(x) = t^2(x-1)(e^{tx} - 1)$. Consequently, $\delta'(x)$ has the same sign as $x(x-1)$, which leads to (3.11). Since the random variables $s^k(X_k)$ are centered, taking $x = s^k(x_k)$ and integrating with respect to the marginal law of $X_k$, we get Lemma 3.3. □

From Lemma 3.3 and (3.10) we have

(3.12) $\quad E_n^k(h_k r(t, g_k/h_k)) \leq \frac{t^2}{2} E_n^k\left(\sum_i \mathbb{1}_{\tau=i} \exp(tS_n^k(s_i)) \mathbb{E}(s_i^k(X_k))^2\right).$

Now $\exp(tS_n^k(s_i)) \leq \exp(2t + tS_n(s_i))$, and therefrom

$$\sum_{k=1}^n \mathbb{E}(h_k r(t, g_k/h_k)) \leq \frac{t^2 e^{2t}}{2}\mathbb{E}\left(\sum_i \mathbb{1}_{\tau=i} \exp(tS_n(s_i)) \sum_{k=1}^n \mathbb{E}(s_i^k(X_k))^2\right).$$

Since $\sum_k \mathbb{E}(s_i^k(X_k))^2 \leq V_n$, we infer that

(3.13) $\qquad\qquad\qquad \sum_{k=1}^n \mathbb{E}(h_k r(t, g_k/h_k)) \leq \tfrac{1}{2}t^2 e^{2t} V_n F(t).$

Together with Lemma 3.2 and (3.8), (3.13) leads to the differential inequality

(3.14) $$tL' - (te^{2t} + 1)L \leq t^2 e^{2t}(V_n/2).$$

Let $\gamma(t) = t^{-2}\exp((1-e^{2t})/2)$. Multiplying (3.14) by $\gamma$, we get

(3.15) $$(t\gamma L)' \leq (V_n/2)e^{2t}\exp((1-e^{2t})/2).$$

Since $t\gamma(t) \sim (1/t)$ as $t$ tends to 0, integrating (3.15) gives

$$t\gamma(t)L(t) \leq \mathbb{E}(Z) + (V_n/2)(1 - \exp((1-e^{2t})/2)),$$

which implies Theorem 1.1(a).

To prove Theorem 1.1(b), we apply both Markov's inequality to the random variable $\exp(tZ)$ and Theorem 1.1(a) with $t = \frac{1}{2}\log(1 + 2\log(1 + x/v))$.

To prove Theorem 1.1(c), we bound up the log-Laplace transform of $Z - \mathbb{E}(Z)$ via Lemma 3.4 and next we apply Markov's exponential inequality.

LEMMA 3.4. *Under the assumptions of Theorem 1.1, for any $t$ in $]0, 2/3[$,*

$$L(t) \leq t\mathbb{E}(Z) + (2\mathbb{E}(Z) + V_n)\frac{t^2}{2 - 3t}.$$

PROOF. From Theorem 1.1(a), it is enough to prove that

$$\exp((e^{2t} - 1)/2) \leq 1 + 2t/(2 - 3t).$$

This inequality holds if and only if

$$\lambda(t) := \log(2 - t) - \log(2 - 3t) - (e^{2t} - 1)/2 \geq 0.$$

Expanding $\lambda$ in power series yields $\lambda(t) = \sum_{j \geq 2} b_j t^j/j!$, where

$$b_j = (j-1)!((3/2)^j - (1/2)^j) - 2^{j-1} \geq 2(j-1)! - 2^{j-1} \geq 0.$$

Hence $\lambda(t) \geq 0$, which implies Lemma 3.4. □

Theorem 1.1(c) follows from Lemma 3.4 by noting that the Legendre transform of the function $t \to t^2/(2 - 3t)$ (here $t < 2/3$) is equal to $\frac{4}{9}(1 + (3x/2) - \sqrt{1 + 3x})$.

**4. Compensated empirical processes.** In this section we prove Theorem 1.2. We start by proving Theorem 1.2(a). Throughout the section, $t$ is any positive real. For $i$ in $\{1, \ldots, m\}$, let

$$L_i(t) = \log \mathbb{E}(\exp(-tS_n(s_i))).$$

Let us define the exponentially compensated empirical process $T(s_i, t)$ by

(4.1) $$T(s_i, t) = S_n(s_i) + t^{-1}L_i(t).$$

We set

(4.2) $$Z_t = \sup_{1 \leq i \leq m} T(s_i, t) \quad \text{and} \quad f_t = \exp(-tZ_t).$$

Let

(4.3) $$F(t) = \mathbb{E}(f_t) = \mathbb{E}(\exp(-tZ_t)) \quad \text{and} \quad \Lambda(t) = \log F(t).$$

Our purpose is to obtain a differential inequality for $\Lambda$ via the log-Sobolev method.

Before that, we link the log-Laplace $L_{-Z}$ of $-Z$ with $\Lambda$.

LEMMA 4.1. *For any positive $t$,*

$$L_{-Z}(t) - \sup_i L_i(t) \leq \Lambda(t) \leq \min(L_{-Z}(t), 0).$$

PROOF. By definition of $Z_t$,

$$\exp(-tZ_t) = \exp\left(\inf_i(-tS_n(s_i) - L_i(t))\right) \geq \exp\left(-tZ - \sup_i L_i\right).$$

Consequently, for any positive $t$,

$$\exp(\Lambda(t)) \geq \exp\left(-\sup_i L_i(t)\right) \mathbb{E} \exp(-tZ),$$

which gives the first inequality. Next, by definition of $Z_t$,

$$\exp(\Lambda(t)) = \mathbb{E}\left(\inf_i \exp(-tS_n(s_i) - L_i(t))\right) \leq \mathbb{E}\left(\exp(-tS_n(s_1) - L_1(t))\right) = 1,$$

which ensures that $\Lambda(t) \leq 0$. Moreover, $L_i(t) \geq 0$ by the centering assumption on the random variables $S_n(s)$. Hence,

$$\exp(\Lambda(t)) \leq \mathbb{E}\left(\inf_i \exp(-tS_n(s_i))\right) = \mathbb{E}(\exp(-tZ)),$$

which completes the proof of Lemma 4.1. □

DEFINITION 4.1. Let $\tau_t$ denote the first integer $i$ such that $Z_t = T(s_i, t)$, where $Z_t$ is defined in (4.2).

Since the random functions $T(s_i, t)$ are analytic functions of $t$, the random function $f_t$ defined in (4.2) is continuous and piecewise analytic, with derivative with respect to $t$, almost everywhere (a.e.):

(4.4) $$f'_t = -Z_t f_t - (L'_{\tau_t}(t) - t^{-1} L_{\tau_t}(t)) f_t = -(Z_t + tZ'_t) f_t,$$

where $tZ'_t = L'_{\tau_t}(t) - t^{-1}L_{\tau_t}(t)$ by convention. Consequently, the Fubini theorem applies and

$$(4.5) \qquad F(t) = 1 - \int_0^t \mathbb{E}((Z_u + uZ'_u)f_u)\,du.$$

Therefrom the function $F$ is absolutely continuous with respect to the Lebesgue measure, with a.e. derivative in the sense of Lebesgue

$$(4.6) \qquad F'(t) = -\mathbb{E}((Z_t + tZ'_t)f_t).$$

Moreover, from the elementary lower bound $f_t \geq \exp(-2nt)$, the function $\Lambda = \log F$ is absolutely continuous with respect to the Lebesgue measure, with a.e. derivative $F'/F$ if $F'$ is the above defined function.

DEFINITION 4.2. Let $f^k = E_n^k f_t$.

We now apply Proposition 2.1 to the random function $f_t$. Clearly,

$$(4.7) \qquad \begin{aligned}\mathbb{E}(f_t \log f_t) - \mathbb{E}(f_t)\log \mathbb{E}(f_t) \\ = \mathbb{E}(t^2 Z'_t f_t) + tF'(t) - F(t)\log F(t) \qquad \text{a.e.}\end{aligned}$$

Hence, applying Proposition 2.1 with $f = f_t$,

$$(4.8) \qquad \begin{aligned}tF' - F\log F \leq -\mathbb{E}(t^2 Z'_t f) + \sum_{k=1}^n \mathbb{E}(g_k \log(g_k/E_n^k g_k)) \\ + \sum_{k=1}^n \mathbb{E}((g_k - f)\log(f^k/f)).\end{aligned}$$

Now choose

$$(4.9) \qquad g_k = \sum_i P_n^k(\tau_t = i)\exp(-tS_n(s_i) - L_i(t)).$$

By definition of $Z_t$,

$$\exp(-tS_n(s_i) - L_i(t)) \geq \exp(-tZ_t),$$

which implies that $g_k \geq f$. Therefore the upper bound on the second term in (4.8) will be derived from Lemma 4.2.

NOTATION 4.1. For sake of brevity, throughout we note $\tau = \tau_t$ and $f_t$.

LEMMA 4.2. Let $\psi(t) = (\exp(2t)+1)/2$. Set $l_{ki}(t) = \log \mathbb{E}(\exp(-ts_i^k(X_k)))$. Then a.s.

$$(f^k/f) \leq \exp(ts_\tau^k(X_k) + l_{k\tau}(t)) \leq \psi(t).$$

PROOF. Let $S_n^k(s) = S_n(s) - s^k(X_k)$. Set

$$Z^k = \sup\{S_n^k(s) + t^{-1} \log \mathbb{E}(\exp(-tS_n^k(s))) : s \in \mathcal{S}\}.$$

Let $\tau_k$ be the first integer in $[1,m]$ such that

$$S_n^k(s_{\tau_k}) + t^{-1} \log \mathbb{E}(\exp(-tS_n^k(s_{\tau_k}))) = Z^k.$$

Clearly

$$f_t \leq \exp(-tZ^k) \exp(-ts_{\tau_k}^k(X_k)) - l_{k\tau_k}(t).$$

Since the stopping time $\tau_k$ is $\mathcal{F}_n^k$-measurable, it follows that

(4.10) $$E_n^k f_t \leq \exp(-tZ^k).$$

Now, by definition of $Z^k$,

(4.11) $$\exp(-tZ^k) \leq \exp(-tZ + ts_\tau^k(X_k) + l_{k\tau}(t)),$$

which ensures that $(f^k/f) \leq \exp(ts_\tau^k(X_k) + l_{k\tau}(t))$. To conclude the proof of Lemma 4.2, recall that $\mathbb{E}(\exp(tX)) \leq \cosh(t)$, for any centered random variable $X$ with values in $[-1,1]$, which implies the second part of Lemma 4.2. □

The next step to bound up the second term on the right-hand side is Lemma 4.3. However, due to technical difficulties, we are able to bound up this term only on some finite interval.

LEMMA 4.3. *Let $(g_k)$ be the finite sequence of random variables defined in* (4.9). *Set $\varphi = \psi \log \psi$. Let $t_0$ be the positive solution of the equation $\varphi(t) = 1$. Then, for any $t$ in $[0, t_0[$,*

$$\sum_{k=1}^n \mathbb{E}((g_k - f) \log(f^k/f)) \leq \frac{\varphi(t)}{1 - \varphi(t)} \left( \sum_{k=1}^n \mathbb{E}(g_k \log(g_k/E_n^k g_k)) - \mathbb{E}(f \log f) \right).$$

PROOF. Since the random variables $S_n^k(s)$ are $\mathcal{F}_n^k$-measurable,

(4.12) $$E_n^k(g_k) = \sum_{i=1}^m P_n^k(\tau = i) \exp(-tS_n^k(s_i) - L_i(t) + l_{ki}(t)).$$

It follows that $E_n^k(g_k) = E_n^k(f \exp(ts_\tau^k(X_k) + l_{k\tau}(t)))$. Hence

(4.13) $$\sum_{k=1}^n \mathbb{E}(g_k - f) = \sum_{k=1}^n \mathbb{E}(f(\exp(ts_\tau^k(X_k) + l_{k\tau}(t)) - 1)).$$

Setting $\eta_k = ts_\tau^k(X_k) + l_{k\tau}(t)$, we have

$$\sum_{k=1}^n \mathbb{E}(g_k - f) = \sum_{k=1}^n \mathbb{E}(f(e^{\eta_k} - 1 - \psi(t)\eta_k)) + \psi(t)\mathbb{E}\left(f\sum_{k=1}^n \eta_k\right)$$

$$= \sum_{k=1}^n \mathbb{E}(f(e^{\eta_k} - 1 - \psi(t)\eta_k)) - \psi(t)\mathbb{E}(f \log f),$$

since $\sum_{k=1}^n \eta_k = -\log f$. Now, for $x$ in $]-\infty, \log \psi(t)]$, the function $x \to e^x - 1 - x\psi(t)$ is nonincreasing. Since $\log \psi(t) \geq \eta_k \geq \log(f^k/f)$ by Lemma 4.2, we infer that

$$\sum_{k=1}^n \mathbb{E}(g_k - f) \leq \sum_{k=1}^n \mathbb{E}(f((f^k/f) - 1 - \psi(t)\log(f^k/f))) - \psi(t)\mathbb{E}(f \log f)$$

$$\leq \psi(t)\left(\sum_{k=1}^n \mathbb{E}(f \log(f/f^k)) - \mathbb{E}(f \log f)\right).$$

Hence, applying (2.3), we obtain

(4.14)
$$\sum_{k=1}^n \mathbb{E}(g_k - f) \leq \psi(t)\left(\sum_{k=1}^n \mathbb{E}(g_k \log(g_k/E_n^k g_k)) + (g_k - f)\log(f^k/f)) - \mathbb{E}(f \log f)\right).$$

Now, from Lemma 4.2 we know that $\log(f^k/f) \leq \log \psi(t)$. Since $(g_k - f) \geq 0$, it follows that

$$\sum_{k=1}^n \mathbb{E}\left((g_k - f)\log\left(\frac{f^k}{f}\right)\right)$$

$$\leq \log \psi(t) \sum_{k=1}^n \mathbb{E}(g_k - f)$$

(4.15)
$$\leq \varphi(t)\left(\sum_{k=1}^n \mathbb{E}\left(g_k \log\left(\frac{g_k}{E_n^k g_k}\right)\right. + (g_k - f)\log\left(\frac{f^k}{f}\right)\right) - \mathbb{E}(f \log f)\right).$$

Since $1 - \varphi(t) > 0$ for any $t$ in $[0, t_0[$, inequality (4.15) then implies Lemma 4.3.

□

From Lemma 4.3 and the differential inequality (4.8) we then get that

$$(1-\varphi)(tF' - F\log F)$$
$$\leq \varphi \mathbb{E}(t^2 Z'_t f - f\log f) - \mathbb{E}(t^2 Z'_t f) + \sum_{k=1}^{n} \mathbb{E}(g_k \log(g_k/E_n^k g_k)),$$

where $\varphi = \varphi(t)$. Now from (4.7), $\mathbb{E}(t^2 Z'_t f - f\log f) = -tF'$, whence

$$(4.16) \quad tF' - (1-\varphi)F\log F \leq -\mathbb{E}(t^2 Z'_t f) + \sum_{k=1}^{n} \mathbb{E}(g_k \log(g_k/E_n^k g_k)).$$

Let us now bound up the first term on the right-hand side in (4.16). Set $w_k = (g_k/E_n^k g_k)$. Then

$$E_n^k(g_k \log(g_k/E_n^k g_k)) = E_n^k(g_k) E_n^k(w_k \log w_k).$$

From (4.12), by convexity of the function $x\log x$,

$$E_n^k(g_k) w_k \log w_k \leq \sum_i P_n^k(\tau = i)(-ts_i^k(X_k) - l_{ki}(t))\exp(-tS_n(s_i) - L_i(t)).$$

Consequently

$$E_n^k(g_k \log(g_k/E_n^k g_k))$$
$$\leq \sum_i P_n^k(\tau = i)\exp(-tS_n^k(s_i) - L_i(t) + l_{ki}(t))(tl'_{ki}(t) - l_{ki}(t)).$$

Since

$$\sum_i P_n^k(\tau = i)\exp(-tS_n^k(s_i) - L_i + l_{ki})(tl'_{ki} - l_{ki})$$
$$= E_n^k\left(\sum_i \mathbb{1}_{\tau=i}\exp(-tS_n^k(s_i) - L_i + l_{ki})(tl'_{ki} - l_{ki})\right),$$

it implies that

$$\mathbb{E}(g_k \log(g_k/E_n^k g_k)) \leq \mathbb{E}(\exp(-tZ_t + ts_\tau^k(X_k) + l_{k\tau})(tl'_{k\tau} - l_{k\tau})).$$

From the convexity of the functions $l_{ki}$, we know that $tl'_{k\tau} - l_{k\tau} \geq 0$. Hence, applying Lemma 4.2, we get

$$\mathbb{E}(g_k \log(g_k/E_n^k g_k)) \leq \psi(t)\mathbb{E}((tl'_{k\tau} - l_{k\tau})f).$$

Since $t^2 Z'_t = tL'_\tau - L_\tau$, it follows that

$$(4.17) \quad -\mathbb{E}(t^2 Z'_t f) + \sum_{k=1}^{n} \mathbb{E}(g_k \log(g_k/E_n^k g_k)) \leq (\psi(t) - 1)\mathbb{E}((tL'_\tau - L_\tau)f).$$

Both (4.16) and (4.17) yield, for $t$ in $[0, t_0[$,

(4.18) $\quad tF' - (1-\varphi)F \log F \leq (\psi(t) - 1)\mathbb{E}((tL'_\tau - L_\tau)f).$

Since $tL'_\tau - L_\tau \leq \sup_i(tL'_i - L_i)$, dividing by $F$, we infer that

(4.19) $\quad t\Lambda' - (1-\varphi)\Lambda \leq (\psi(t) - 1)\sup_i(tL'_i - L_i).$

Next we derive an upper bound on $tL'_i - L_i$ from Lemma 4.4.

LEMMA 4.4. *Let $Y$ be a random variable with values in $]-\infty, 1]$, such that $\mathbb{E}(Y^2) < +\infty$. Then, for any positive $t$,*

$$\mathbb{E}(tYe^{tY}) - \mathbb{E}(e^{tY})\log \mathbb{E}(e^{tY}) \leq \mathbb{E}(Y^2)(1 + (t-1)e^t).$$

PROOF. From the variational definition of the entropy in Ledoux (1996) we know that, for any positive constant $c$ and any positive random variable $T$,

$$\mathbb{E}(T\log T) - \mathbb{E}(T)\log \mathbb{E}(T) \leq \mathbb{E}(T\log(T/c) - T + c).$$

Taking $c = 1$ and $T = \exp(tY)$, we then get that

$$\mathbb{E}(tYe^{tY}) - \mathbb{E}(e^{tY})\log \mathbb{E}(e^{tY}) \leq \mathbb{E}((tY - 1)e^{tY} + 1).$$

Now, from l'Hôpital's rule for monotonicity the function $x \to x^{-2}(1 + (x-1)e^x)$ is nondecreasing on the real line. Hence, for any positive $t$,

$$(tY - 1)e^{tY} + 1 \leq Y^2(1 + (t-1)e^t),$$

which implies Lemma 4.4. □

Let $Y_k = s_i^k(X_k)$. From the centering assumption, $\mathbb{E}\exp(tY_k) \geq 1$. Hence we have

$$tl'_{ki}(t) - l_{ki}(t) \leq \mathbb{E}(tY_k e^{tY_k}) - \mathbb{E}(e^{tY_k})\log \mathbb{E}(e^{tY_k}) \leq (1 + (t-1)e^t)\operatorname{Var} Y_k$$

by Lemma 4.4. Since $tL_i - l'_i = \sum_k(tl'_{ki} - l_{ki})$, it ensures that

(4.20) $\quad tL'_i - L_i \leq V_n(1 + (t-1)e^t).$

Both the above bound and (4.19) lead to the differential inequality below.

PROPOSITION 4.1. *For any $t$ in $[0, t_0[$,*

(4.21) $\quad t\Lambda' - (1-\varphi)\Lambda \leq \tfrac{1}{2}V_n(e^{2t} - 1)(1 + (t-1)e^t).$

It remains to bound up $\Lambda$. Set $\widetilde{\Lambda}(t) = t^{-1}\Lambda(t)$ and

$$(4.22) \qquad I(t) = \int_0^t \frac{\varphi(u)}{u}\, du.$$

Then $(\widetilde{\Lambda} e^I)' = t^{-2}(t\Lambda' - (1-\varphi)\Lambda)e^I$. Consequently, from Proposition 4.1,

$$(4.23) \qquad (\widetilde{\Lambda} e^I)' \leq \frac{V_n}{2t^2}(e^{2t} - 1)(1 + (t-1)e^t)e^I.$$

Since $\widetilde{\Lambda} e^I$ is absolutely continuous with respect to the Lebesgue measure, integrating (4.23) yields

$$(4.24) \qquad \begin{aligned} \widetilde{\Lambda}(t) &\leq \widetilde{\Lambda}(\varepsilon) e^{I(\varepsilon) - I(t)} \\ &\quad + \frac{V_n}{2}\int_\varepsilon^t u^{-2}(e^{2u} - 1)(1 + (u-1)e^u)e^{I(u) - I(t)}\, du \end{aligned}$$

for $0 < \varepsilon < t$. The control of the integral on the right-hand side will be done via the bounds for $\varphi$ below, whose proof is carried out in Section 5.

LEMMA 4.5. *For any $t$ in $[0, t_0]$, $t \leq \varphi(t) \leq t\exp(2t) - (t^2/2)$.*

By Lemma 4.5, $\lim_0 I(\varepsilon) = 0$. Furthermore, $\widetilde{\Lambda}(\varepsilon) \leq \varepsilon^{-1} L_{-Z}(\varepsilon)$ by Lemma 4.1. Therefore

$$(4.25) \qquad \limsup_{\varepsilon \to 0} \widetilde{\Lambda}(\varepsilon) e^{I(\varepsilon) - I(t)} \leq -\mathbb{E}(Z) e^{-I(t)}.$$

Now $I(u) - I(t) \leq (u - t)$ by Lemma 4.5. Consequently, letting $\varepsilon \to 0$ in (4.24) and applying (4.25), we get

$$(4.26) \quad \Lambda(t) \leq -\mathbb{E}(Z) t e^{-I(t)} + \tfrac{1}{2} V_n t e^{-t} \int_0^t u^{-2}(e^{2u} - 1)(1 + (u-1)e^u)e^u\, du.$$

To bound up $L_{-Z}$, we then apply the Bennett bound $L_i(t) \leq V_n(e^t - t - 1)$ together with Lemma 4.1. This yields the Proposition 4.2.

PROPOSITION 4.2. *Let the function $J$ be defined by*

$$J(t) = \tfrac{1}{2}\int_0^t u^{-2}(e^{2u} - 1)(1 + (u-1)e^u)e^u\, du$$

*and let $I$ be the function defined in (4.22). For any $t$ in $[0, t_0]$,*

$$L_{-Z}(t) + t\mathbb{E}(Z) \leq t\mathbb{E}(Z)(1 - e^{-I(t)}) + V_n(t e^{-t} J(t) + e^t - t - 1).$$

To obtain Theorem 1.2(a) for $t$ in $[0, t_0]$, we bound up the functions appearing in Proposition 4.2 via Lemma 4.6, proved in Section 5.

LEMMA 4.6. *For any $t$ in $[0, t_0]$,*

(a) $$te^{-t}J(t) + e^t - t - 1 \leq \tfrac{1}{9}(\exp(3t) - 3t - 1),$$

(b) $$t(1 - e^{-I(t)}) \leq \tfrac{2}{9}(\exp(3t) - 3t - 1).$$

Next, proceeding as in Klein (2002), we prove (a) for $t$ in $[t_0, +\infty[$. For sake of brevity, set $E = \mathbb{E}(Z)$. By Lemma 4.1, for any positive $t$,

$$L_{-Z}(t) + tE \leq tE + \sup_i L_i(t) \leq tE + V_n(e^t - t - 1) \leq v\max(t/2, e^t - t - 1).$$

Now, let $t_1$ be the unique positive solution of the equation $e^t - t - 1 = t/2$. $t_1$ belongs to $[0.76, 0.77]$, whence $t_1 > t_0$ (note that $t_0 \in [0.46, 0.47]$). If $t \geq t_1$, then $t/2 \leq e^t - t - 1$. In that case

$$L_{-Z}(t) + tE \leq v(e^t - t - 1) \leq (v/9)(e^{3t} - 3t - 1),$$

which proves (a) for $t \geq t_1$.

If $t$ belongs to $[t_0, t_1]$, from the convexity of $L_{-Z}$ we have

$$L_{-Z}(t) + tE \leq \frac{1}{9}v(e^{3t_0} - 3t_0 - 1)\frac{t_1 - t}{t_1 - t_0} + \frac{1}{2}vt_1\frac{t - t_0}{t_1 - t_0} \leq \frac{1}{9}v(e^{3t} - 3t - 1),$$

which completes the proof of Theorem 1.2(a).

To prove Theorem 1.2(c), we note that, for any $t$ in $[0, 1[$,

(4.27) $$\frac{1}{9}(e^{3t} - 3t - 1) \leq \frac{t^2}{2 - 2t}$$

[cf. Rio (2000), page 152]. Theorem 1.2(c) follows from (4.27) via the usual Cramér–Chernoff calculation.

**5. Technical tools.** In this section, we prove Lemmas 4.5 and 4.6.

PROOF OF LEMMA 4.5. By definition of $\psi$ and $\varphi$,

$$\varphi(t) = t\psi(t) + \psi(t)\log\cosh(t) \geq t,$$

since $\psi(t) \geq 1$ for any nonnegative $t$. Next

$$t\exp(2t) - (t^2/2) - \varphi(t) = \psi(t)(t\tanh(t) - \log\cosh(t) - (e^{2t} + 1)^{-1}t^2),$$

so that Lemma 4.5 holds if $p(t) := t\tanh(t) - \log\cosh(t) - t^2/(e^{2t} + 1) \geq 0$ for $t$ in $[0, t_0]$. Now $p(0) = 0$ and

$$2\cosh^2(t)p'(t) = 2t - t(e^{-2t} + 1) - t^2 = t(1 - t - e^{-2t}).$$

Since $\exp(-2t) \leq 1 - t$ for $t$ in $[0, 1/2]$, the above identity ensures that $p'(t) \geq 0$ on $[0, t_0]$ (recall that $t_0 < 1/2$), which implies Lemma 4.5. □

PROOF OF LEMMA 4.6. We start by proving (a). Clearly (a) holds if

(5.1) $\quad \alpha(t) = \frac{1}{9}e^t(e^{3t} - 3t - 1) + e^t(1 + t - e^t) - J(t) \geq 0$

for any $t$ in $[0, 4]$ [with the convention $\alpha(0) = 0$]. The function $\alpha$ is analytic on the real line. To prove (5.1), we then note that $\alpha^{(i)}(0) = 0$ for $i = 1, 2$. Consequently (a) holds if, for $t$ in $[0, 4]$,

(5.2) $\quad \alpha^{(3)}(t) = e^{3t}(-t + (19/3)) - 4e^{2t} + e^t(-3t - 5) + (8/3) > 0.$

Now (5.2) holds if $\alpha^{(4)}(t) > 0$, since $\alpha^{(3)}(0) > 0$. Next

$$\beta(t) := e^{-t}\alpha^{(4)}(t) = 3e^{2t}(-2t + 11) - 8e^t - 3$$

satisfies $\beta(0) > 0$ and, for $t$ in $[0, 4]$,

$$\beta'(t) = 12e^{2t}(5 - t) - 8e^t > e^t(12e^t - 8) > 0,$$

which ensures that $\beta(t) > 0$ for $t$ in $[0, 4]$. Hence Lemma 4.6(a) holds.

To prove (b), we apply Lemma 4.5 to bound up the function $I(t)$. This gives

$$I(t) \leq \int_0^t (e^{2u} - u/2)\, du = \frac{e^{2t} - 1}{2} - \frac{t^2}{4}.$$

Now, recall $t_0 \leq 1/2$. For $t$ in $[0, 1/2]$, expanding $\exp(2t)$ in entire series yields

$$(\exp(2t) - 1)/2 = t + t^2 + 4t^3 \sum_{k \geq 3} \frac{1}{k!}(2t)^{k-3} \leq t + t^2 + 4t^3 \sum_{k \geq 3} \frac{1}{k!}.$$

Hence, for $t \leq 1/2$,

(5.3) $\quad I(t) \leq t + \frac{3}{4}t^2 + (4e - 10)t^3 \leq t + \frac{3}{4}t^2 + \frac{7}{8}t^3 =: \gamma(t).$

From (5.3), Lemma 4.6(b) holds if

$$d(t) = \frac{2}{9}(e^{3t} - 3t - 1) - t + t\exp(-\gamma(t)) \geq 0.$$

Now $d(0) = d'(0) = 0$ and

$$d''(t) = 2e^{-\gamma(t)}(e^{4t + (3/4)t^2 + (7/8)t^3} - 1 - \tfrac{7}{4}t - \tfrac{15}{4}t^2 + \tfrac{15}{4}t^3 + \tfrac{63}{16}t^4 + \tfrac{441}{128}t^5).$$

Since

$$e^{4t + (3/4)t^2 + (7/8)t^3} \geq e^{4t} \geq 1 + 4t + 8t^2,$$

we have $d''(t) > 0$ for any positive $t$. Consequently, $d(t) \geq 0$, which implies Lemma 4.6(b). □

## REFERENCES

Bobkov, S. (1996). Some extremal properties of the Bernoulli distribution. *Theory Probab. Appl.* **41** 748–755. MR1687168

Bousquet, O. (2003). Concentration inequalities for sub-additive functions using the entropy method. *Stochastic Inequalities and Applications* **56** 213–247. MR2073435

Klein, T. (2002). Une inégalité de concentration à gauche pour les processus empiriques. *C. R. Acad. Sci. Paris Sér. I Math.* **334** 495–500. MR1890640

Ledoux, M. (1996). On Talagrand's deviation inequalities for product measures. *ESAIM Probab. Statist.* **1** 63–87. MR1399224

Ledoux, M. and Talagrand, M. (1991). *Probability in Banach Spaces. Isoperimetry and Processes*. Springer, Berlin. MR1102015

Massart, P. (2000). About the constants in Talagrand's concentration inequalities for empirical processes. *Ann. Probab.* **28** 863–884. MR1782276

Panchenko, D. (2001). A note on Talagrand's concentration inequality. *Electron. Comm. Probab.* **6** 55–65. MR1831801

Panchenko, D. (2003). Symmetrization approach to concentration inequalities for empirical processes. *Ann. Probab.* **31** 2068–2081. MR2016612

Rio, E. (2000). Théorie asymptotique des processus aléatoires faiblement dépendants. In *Mathématiques et Applications* (J. M. Ghidaglia et X. Guyon, eds.) **31**. Springer, Berlin. MR2117923

Rio, E. (2001). Inégalités de concentration pour les processus empiriques de classes de parties. *Probab. Theory Related Fields* **119** 163–175. MR1818244

Rio, E. (2002). Une inégalité de Bennett pour les maxima de processus empiriques. *Ann. Inst. H. Poincaré Probab. Statist.* **38** 1053–1057. MR1955352

Talagrand, M. (1996). New concentration inequalities in product spaces. *Invent. Math.* **126** 503–563. MR1419006

UMR 8100 CNRS
Université de Versailles
 Saint Quentin en Yvelines
45 Avenue des Etats Unis
78035 Versailles
France
e-mail: rio@math.uvsq.fr
e-mail: klein@math.uvsq.fr